\documentclass[12pt]{article}
\usepackage{amssymb}
\usepackage{latexsym,bm}
\usepackage{graphicx}
\usepackage{amsmath}
\usepackage{mathrsfs}
\usepackage{epstopdf}

\date{}
\begin{document}
\title{A sufficient condition for pancyclic graphs}
\author{\hskip -10mm Xingzhi Zhan\footnote{E-mail address: \tt zhan@math.ecnu.edu.cn}\\
{\hskip -10mm \small Department of Mathematics, East China Normal University, Shanghai 200241, China}}\maketitle
\begin{abstract}
A graph $G$ is called an $[s,t]$-graph if any induced subgraph of $G$ of order $s$ has size at least $t.$ We prove that every $2$-connected $[4,2]$-graph of order
at least $7$ is pancyclic. This strengthens existing results. There are $2$-connected $[4,2]$-graphs which do not satisfy the Chv\'{a}tal-Erd\H{o}s condition.
We also determine the triangle-free graphs among $[p+2,p]$-graphs for a general $p.$
\end{abstract}

{\bf Key words.} Hamiltonian graph; pancyclic graph; $[s,t]$-graph; triangle-free

{\bf Mathematics Subject Classification.} 05C38, 05C42, 05C45, 05C75
\vskip 8mm

\section{Introduction}

We consider finite simple graphs and use standard terminology and notation from [3] and [8]. The {\it order} of a graph is its number of vertices, and the
{\it size} its number of edges.  A $k$-cycle is a cycle of length $k.$ In 1971 Bondy [1] introduced the concept of a pancyclic graph. A graph $G$ of order $n$ is called
{\it pancyclic} if for every integer $k$ with $3\le k\le n,$ $G$ contains a $k$-cycle.

{\bf Definition 1.} Let $s$ and $t$ be given integers. A graph $G$ is called an {\it $[s,t]$-graph} if any induced subgraph of $G$ of order $s$ has size at least $t.$

Denote by $\alpha (G)$ the independence number of a graph $G$. We have two facts. (1) Every $[s,t]$-graph is an $[s+1,t+1]$-graph; (2) $\alpha(G)\le k$ if and only if
$G$ is a $[k+1,1]$-graph. Thus the concept of an $[s,t]$-graph is an extension of the independence number.

In 2005 Liu and Wang [6] proved the following result.

{\bf Theorem 1.} {\it Every $2$-connected $[4,2]$-graph of order at least $6$ is hamiltonian.}

In 2007 Liu, Wang and Gao [7] improved Theorem 1 as follows.

{\bf Theorem 2.} {\it Let $G$ be a $2$-connected $[4,2]$-graph of order $n$ with $n\ge 7.$ If $G$ contains a $k$-cycle with $k<n,$ then $G$ contains a $(k+1)$-cycle.}

In this paper we further strengthen Theorem 2 by proving that every $2$-connected $[4,2]$-graph of order at least $7$ is pancyclic (Theorem 6). To do so, we will determine
the triangle-free graphs among $[p+2,p]$-graphs. This preliminary result (Lemma 5) is of independent interest.

\section{Main results}

We denote by $V(G)$ and $E(G)$ the vertex set and edge set of a graph $G,$ respectively, and denote by $|G|$ and $e(G)$ the order and size of $G,$ respectively.
Thus $|G|=|V(G)|$ and $e(G)=|E(G)|.$ For a vertex subset $S\subseteq V(G),$ we use $G[S]$ to denote the subgraph of $G$ induced by $S.$ The neighborhood of a vertex
$x$ is denoted by $N(x)$ and the closed neighborhood of $x$ is $N[x]\triangleq N(x)\cup \{x\}.$ The degree of $x$ is denoted by ${\rm deg}(x).$
For $S\subseteq V(G),$ $N_S(x)\triangleq N(x)\cap S$ and the degree of $x$ in $S$ is ${\rm deg}_S(x)\triangleq |N_S(x)|.$ Given two vertex subsets $S$ and $T$ of $G,$
we denote by $[S, T]$ the set of edges having one endpoint in $S$ and the other in $T.$ The degree of $S$ is ${\rm deg}(S)\triangleq |[S, \overline{S}]|,$ where $\overline{S}=V(G)\setminus S.$ We denote by $C_n$ and $K_n$  the cycle of order $n$ and the complete graph of order $n,$ respectively. $\overline{G}$ denotes
the complement of a graph $G.$

We will need the following two lemmas on integral quadratic forms.

{\bf Lemma 3.} {\it Given positive integers $n\ge k\ge 2,$ let $x_1,x_2,\ldots,x_k$ be positive integers such that $\sum_{i=1}^k x_i =n.$ Then
$$
 n-1\le \sum_{i=1}^{k-1}x_ix_{i+1}\le\begin{cases}\lfloor n/2\rfloor\cdot \lceil n/2\rceil\quad {\rm if}\,\,\,k=2,\,3\\
                    ab+k-5 \quad\quad\,\, {\rm if}\,\,\,k\ge 4 \end{cases}\eqno (1)
$$
where $a=\lfloor (n-k+4)/2\rfloor$ and $b=\lceil (n-k+4)/2\rceil.$ For any $n$ and $k,$ the lower and upper bounds in (1) can be attained.
}

{\bf Proof.} Define a quadratic polynomial $f(x_1,x_2,\ldots,x_k)=\sum_{i=1}^{k-1}x_ix_{i+1}.$ We first prove the left-hand side inequality in (1).
Let $x_j={\rm min}\{x_i \,|\, 1\le i\le k\}.$ We have
\begin{align*}
f(x_1,x_2,\ldots,x_k)&\ge x_1x_j+\cdots+x_{j-1}x_j+x_jx_{j+1}+x_jx_{j+2}+\cdots+x_jx_k\\
                     &=x_j(x_1+\cdots+x_{j-1}+x_{j+1}+\cdots+x_k)\\
                     &=x_j(n-x_j)\\
                     &\ge n-1.
\end{align*}
This proves the first inequality in (1). The lower bound $n-1$ is attained at $x_1=n-k+1,$ $x_2=\cdots=x_k=1.$

Now we prove the second inequality in (1). The case $k=2$ is an elementary fact: $f(x_1,x_2)=x_1 x_2\le \lfloor n/2\rfloor\cdot \lceil n/2\rceil$
where equality holds when $x_1=\lfloor n/2\rfloor$ and $x_2=\lceil n/2\rceil.$ The case $k=3$ reduces to the case $k=2$ as follows:
$$
f(x_1,x_2,x_3)=x_2(x_1+x_3)\le \lfloor n/2\rfloor\cdot \lceil n/2\rceil
$$
where equality holds when $x_2=\lfloor n/2\rfloor$ and $x_1+x_3=\lceil n/2\rceil.$ Next suppose $k\ge 4.$ Denote by $f_{\rm max}$ the maximum value of $f.$
If $x_1>1,$ with $x_1^{\prime}=1,$ $x_2^{\prime}=x_2,$ $x_3^{\prime}=x_3+x_1-1$ and $x_i^{\prime}=x_i$ for $i\ge 4$ we have
$$
f(x_1^{\prime}, x_2^{\prime},\ldots,x_k^{\prime})-f(x_1,x_2,\ldots,x_k)=(x_1-1)x_4>0.
$$
Similarly analyzing the variable $x_k,$ we deduce that $f_{\rm max}$ can only be attained at some $x_1,\ldots,x_k$ with $x_1=x_k=1,$ which we assume now.
With $x_2^{\prime}=1,$ $x_3^{\prime}=x_3,$ $x_4^{\prime}=x_4+x_2-1,$ and $x_i^{\prime}=x_i$ for $i=5,\ldots,k-1$ we have
$$
f(1,1,x_3^{\prime}, x_4^{\prime},\ldots,x_{k-1}^{\prime},1)-f(1,x_2,x_3,\ldots,x_{k-1},1)=(x_2-1)(x_5-1)\ge 0.
$$
Hence $f_{\rm max}$ can be attained at a certain $(1,1,x_3,\ldots,x_{k-1},1).$ Successively applying this argument we deduce that $f_{\rm max}$ can be attained at
$(1,1,\ldots,1,x_{k-2},x_{k-1},1).$ Now $(x_{k-2}+1)+(x_{k-1}+1)=n-k+4.$ We have
\begin{align*}
f(1,1,\ldots,1,x_{k-2},x_{k-1},1)&=(x_{k-2}+1)(x_{k-1}+1)+k-5\\
                                 &\le \lfloor (n-k+4)/2\rfloor\cdot \lceil (n-k+4)/2\rceil+k-5.
\end{align*}
This proves the second inequality in (1). The upper bound is attained at $x_1=x_2=\cdots=x_{k-3}=x_k=1,$ $x_{k-2}=\lfloor (n-k+2)/2\rfloor$ and
$x_{k-1}=\lceil (n-k+2)/2\rceil.$ $\hfill\Box$

{\bf Lemma 4.} [9, Theorem 1] {\it Given positive integers $n\ge k\ge 2,$ let $x_1,x_2,\ldots,x_k$ be positive integers such that $\sum_{i=1}^k x_i =n.$ Then
$$
2n-k\le \sum_{i=1}^k x_ix_{i+1} \eqno(2)
$$
where $x_{k+1}\triangleq x_1.$ For any $n$ and $k,$ the lower bound in (2) can be attained.}

The sharp upper bound on the quadratic form in (2) is also determined in [9], but we do not need it here.

{\bf Definition 2.} Given a graph $H$ and a positive integer $k,$ the {\it $k$-blow-up of $H,$} denoted by $H^{(k)},$ is the graph obtained by replacing
every vertex of $H$ with $k$ different vertices where a copy of $u$ is adjacent to a copy of $v$ in the blow-up graph if and only if $u$ is adjacent to $v$ in $H.$

For example, $C_5^{(2)}$ is depicted in Figure 1.
\begin{figure}[h]
\centering
\includegraphics[width=0.25\textwidth]{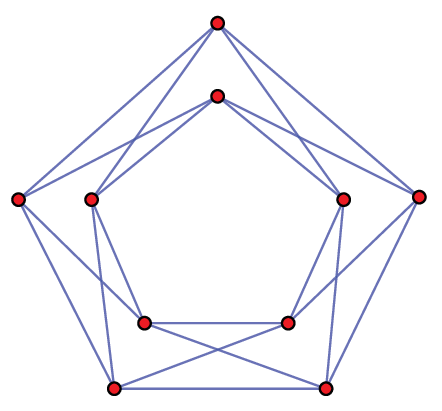}
\caption{The 2-blow-up of $C_5$}
\end{figure}

Now we are ready to determine the triangle-free graphs among $[p+2,p]$-graphs. $\delta(G)$ and $\Delta(G)$ denote the minimum and maximum degrees of a graph $G,$
respectively. We regard isomorphic graphs as the same graph. Thus for two graphs $G$ and $H,$ the notation $G=H$ means that $G$ and $H$ are isomorphic.

{\bf Lemma 5.} {\it Let $G$ be a $[p+2,p]$-graph of order $n$ with $\delta(G)\ge p\ge 2$ and $n\ge 2p+3.$ Then $G$ is triangle-free if and only if
$p$ is even, $p\ge 6$ and $G=C_5^{(p/2)}.$}

{\bf Proof.} We will repeatedly use the condition that $G$ is a $[p+2,p]$-graph without mentioning it possibly.
Denote $\Delta=\Delta(G)$ and choose a vertex $x\in V(G)$ such that ${\rm deg}(x)=\Delta.$ Let $S=N(x)$ and $T=V(G)\setminus S.$ Then $|S|=\Delta.$

Next suppose that $G$ is triangle-free. Then $S$ is an independent set. Since $G$ is a $[p+2,p]$-graph, $\Delta\le p+1.$ We assert that $\Delta=p$ and hence
$G$ is $p$-regular, since $\delta(G)\ge p$ by the assumption. Otherwise $\Delta=p+1.$ Since $n\ge 2p+3,$ $|T|\ge p+2.$ Thus $G[T]$ contains an edge $uv.$
$|\{u\}\cup S|=p+2$ implies that ${\rm deg}_S(u)\ge p.$ Similarly ${\rm deg}_S(v)\ge p.$ Since $p+p=2p>p+1=|S|,$ we have $N_S(u)\cap N_S(v)\neq\emptyset.$
Let $w\in N_S(u)\cap N_S(v).$ Then $wuvw$ is a triangle, a contradiction. This shows that $G$ is $p$-regular.

Let $y\in S$ and denote $C=N(y).$ Then $C$ is an independent set and $|C|=p.$ Denote $D=T\setminus C.$ See the illustration in Figure 2.

\begin{figure}[h]
\centering
\includegraphics[width=0.39\textwidth]{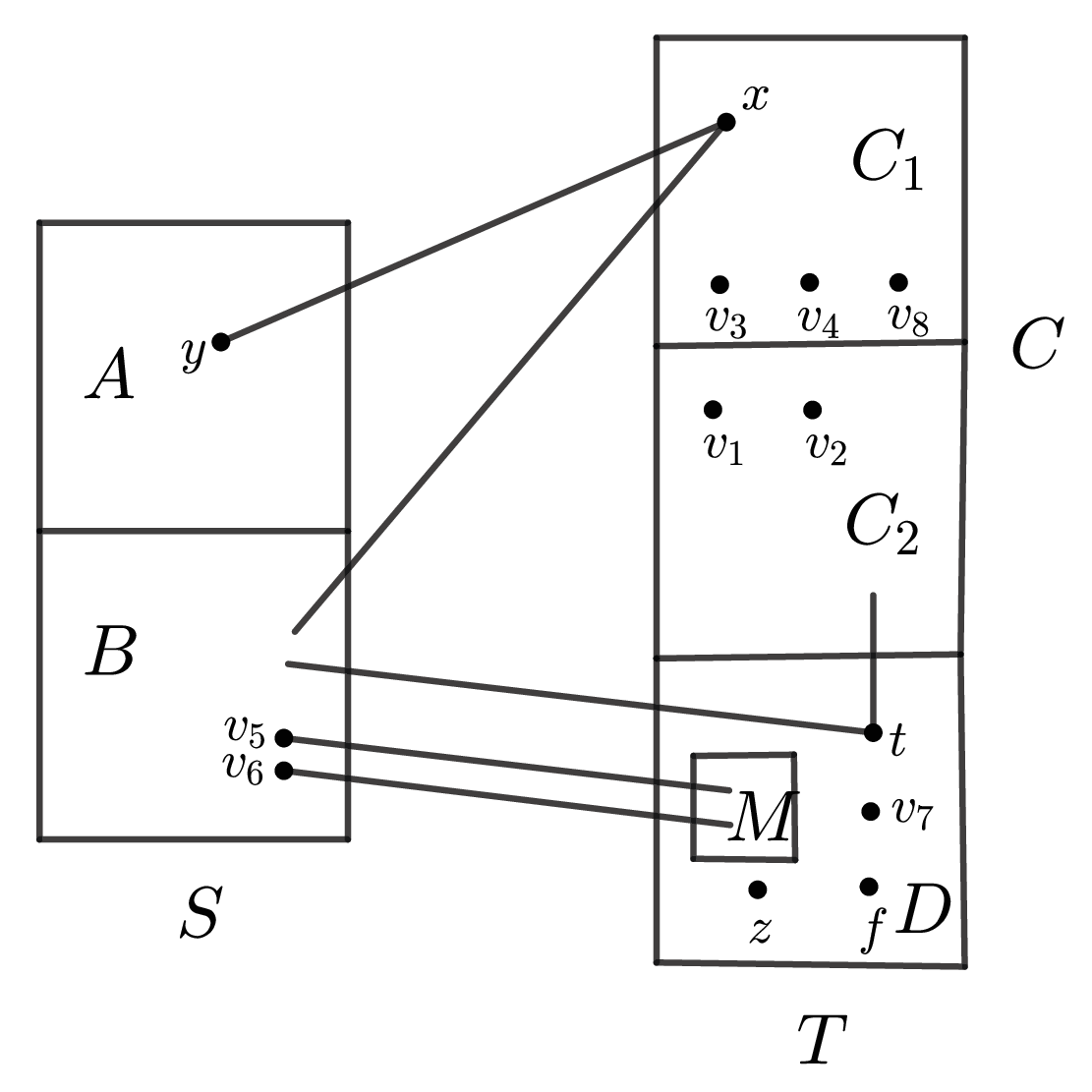}
\caption{The structure of $G$}
\end{figure}

Since $n\ge 2p+3,$ we have $|D|=n-2p\ge 3.$ Thus $D$ is not a clique, since $G$ is triangle-free.
Let $z$ and $f$ be any two distinct nonadjacent vertices in $D.$ Since $|\{z,f\}\cup S|=p+2,$ $|\{z,f\}\cup C|=p+2,$
and $S$ and $C$ are independent sets, we have
$$
{\rm deg}_S(z)+{\rm deg}_S(f)=|[\{z,f\}, S]|\ge p \quad {\rm and} \quad {\rm deg}_C(z)+{\rm deg}_C(f)=|[\{z,f\}, C]|\ge p.
$$
Note that $S\cap C=\emptyset,$ ${\rm deg}(z)={\rm deg}(f)=p.$ We must have
$$
|[\{z,f\}, S]|=p \quad {\rm and} \quad |[\{z,f\}, C]|=p. \eqno(3)
$$
We assert that $D$ is an independent set. Otherwise $D$ contains two adjacent vertices $u_1$ and $u_2.$ Let $u_3\in D\setminus \{u_1, u_2\}.$
Since $G$ is triangle-free, $u_3$ is nonadjacent to at least one vertex in $\{u_1, u_2\},$ say, $u_1$. Setting $z=u_1$ and $f=u_3$ in (3) we deduce that
$|[\{u_1,u_3\}, S\cup C]|=2p.$ On the other hand, since both $u_1$ and $u_3$ have degree $p,$ and $u_1$ already has a neighbor $u_2\not\in S\cup C,$
we have $|[\{u_1,u_3\}, S\cup C]|\le 2p-1,$ a contradiction.

Observe that now (3) holds for any two distinct vertices $z$ and $f$ in $D.$ (3) has the equivalent form
$$
{\rm deg}_S(z)+{\rm deg}_S(f)=p \quad {\rm and} \quad {\rm deg}_C(z)+{\rm deg}_C(f)=p. \eqno(4)
$$
Then (4) and $|D|\ge 3$ imply that for any vertex $z\in D,$
$$
{\rm deg}_S(z)={\rm deg}_C(z)=p/2. \eqno (5)
$$
To see this, to the contrary, first suppose ${\rm deg}_S(z)>p/2.$ Then by the first equality in (4), for any two other vertices $f,r\in D$ we have
${\rm deg}_S(f)<p/2$ and ${\rm deg}_S(r)<p/2,$ yielding ${\rm deg}_S(f)+{\rm deg}_S(r)<p,$ which contradicts (4).
If ${\rm deg}_S(z)<p/2,$ the same argument gives a contradiction. A similar analysis with $C$ in place of $S$ shows ${\rm deg}_C(z)=p/2.$
Thus we have proved (5).
In particular, $q\triangleq p/2$ is a positive integer; i.e., $p$ is even. Now choose an arbitrary but fixed vertex $t\in D$ and denote
$B=N_S(t),$ $C_2=N_C(t),$ $A=S\setminus B$ and $C_1=C\setminus C_2.$ See the illustration in Figure 2. We have
$$
|A|=|B|=|C_1|=|C_2|=q.
$$
Since $G$ is $p$-regular of order $n\ge 2p+3,$ it is impossible that $p=2.$ Otherwise $G$ would be a $2$-regular graph of order $\ge 7,$ which is not a $[4,2]$-graph.
Thus $p\ge 4$ and $q\ge 2.$

Choose any two distinct vertices $v_1, v_2\in C_2.$ $|\{v_1, v_2\}\cup S|=p+2$ implies that $|[\{v_1, v_2\}, S]|\ge p.$ Since $G$ is triangle-free, $N(v_i)\cap B=\emptyset,$
$i=1, 2.$ Hence $N_S(v_i)=N_A(v_i),$ $i=1, 2.$ However, $|A|=q.$ We have $N_S(v_i)=A$ and ${\rm deg}_A(v_i)=q,$ $i=1,2,$ implying that every vertex in $C_2$ is adjacent to
every vertex in $A.$

Choose any two distinct vertices $v_3, v_4\in C_1.$ Then $|\{v_3, v_4\}\cup C_2\cup B|=p+2.$ Since $\{v_3, v_4\}\cup C_2$ is an independent set and
$[C_2, B]=\emptyset,$ we have $|[\{v_3, v_4\}, B]|\ge p.$ However, $|B|=q.$ Hence $N_B(v_j)=B,$ $j=3,4.$ This shows that every vertex in $C_1$ is adjacent to every vertex in $B.$ Consequently, every vertex in $B$ has exactly $q$ neighbors in $D.$

Choose any vertex $v_5\in B.$ Denote $M=N_D(v_5).$ We have $|M|=q.$ Since $G$ is triangle-free and every vertex in $B$ is adjacent to every vertex in $C_1,$ the neighborhood
of any vertex in $M$ is disjoint from $C_1.$ Thus the $q$ neighbors of any vertex of $M$ in $C$ are exactly the vertices of $C_2,$ implying that every vertex in $M$ is adjacent to every vertex in $C_2.$ The neighborhood of any vertex in $C_2$ is $A\cup M.$ For the same reason, for any vertex $v_6\in B$ with $v_6\neq v_5,$ we must have $N_D(v_6)=M.$
Hence the neighborhood of any vertex in $M$ is $B\cup C_2.$

We assert that $M=D.$ Otherwise let $v_7\in D\setminus M.$ Take a vertex $v_8\in C_1.$ Note that $B\cup C_2$ is an independent set of cardinality $p$ and
$[v_7, B\cup C_2]=\emptyset.$ Denote $R=\{v_7, v_8\}\cup B\cup C_2.$ Then $|R|=p+2$ and hence $G[R]$ has size at least $p.$ However, the size of $G[R]$ is at most
$|[v_8, B]|+1=q+1<p,$ a contradiction. Finally, since $G$ is $p$-regular, every vertex in $A$ must be adjacent to every vertex in $C_1.$ Denote
$V_1=A,$ $V_2=C_1,$ $V_3=B,$ $V_4=D,$ $V_5=C_2$ and set $V_6=V_1.$ Then each $V_i$ is an independent set of cardinality $q=p/2$ and every vertex in $V_i$ is adjacent to
every vertex in $V_{i+1}$ for $i=1,2,\ldots,5.$ This proves that $G=C_5^{(q)}.$ Note that we have shown above that $q=|D|\ge 3,$ implying that $p=2q\ge 6.$

Conversely let $H=C_5^{(q)}$ where $q=p/2$ and $p\ge 6$ is even. We will prove that $H$ is a triangle-free $[p+2,p]$-graph. Write $H=H_1\vee H_2\vee H_3\vee H_4\vee H_5\vee H_1$
where each $H_i=\overline{K_q}$ and $\vee$ is the join operation on two vertex-disjoint graphs. If $H$ contains a triangle, it must lie in $H[V(H_i)\cup V(H_{i+1})]$ for some
$i$ ($H_6\triangleq H_1$). However, this is a bipartite graph, containing no triangle.

Let $U\subseteq V(H)$ with $|U|=p+2.$ We need to show $e(H[U])\ge p.$ Denote $I=\{i\, |\, U\cap V(H_i)\neq\emptyset,\,\,\,1\le i\le 5\}.$
Since $|H_i|=q,$ $1\le i\le 5$ and $|U|=p+2,$ we have $|I|\ge 3.$ Denote $x_i=|U\cap V(H_i)|$ for $1\le i\le 5.$ Then $0\le x_i\le q.$ We distinguish three cases.

Case 1. $|I|=3.$

There are at least two consecutive integers in $I$ ($1$ and $5$ are regarded as consecutive here). Without loss of generality, suppose $1,2\in I.$  Then $1\le x_1,\,x_2\le q$ and $x_1+x_2\ge p+2-q=q+2.$
Hence $e(H[U])\ge x_1x_2\ge 2q=p.$

Case 2. $|I|=4.$

Without loss of generality, suppose $I=\{1,2,3,4\}.$ Then $e(H[U])=x_1x_2+x_2x_3+x_3x_4$ where each $x_i$ is a positive integer and $x_1+x_2+x_3+x_4=p+2.$
Applying Lemma 3 we have $e(H[U])\ge (p+2)-1=p+1>p.$

Case 3. $|I|=5.$

Now $e(H[U])=x_1x_2+x_2x_3+x_3x_4+x_4x_5+x_5x_1$ where each $x_i$ is a positive integer and $x_1+x_2+x_3+x_4+x_5=p+2.$ Applying Lemma 4 we have
$e(H[U])\ge 2(p+2)-5=2p-1>p.$

In every case $H$ is a $[p+2,p]$-graph. This completes the proof. $\hfill\Box$

{\bf Theorem 6.} {\it Every $2$-connected $[4,2]$-graph of order at least $7$ is pancyclic.}

{\bf Proof.} Let $G$ be a $2$-connected $[4,2]$-graph of order at least $7.$ Since $G$ is $2$-connected, $\delta (G)\ge 2.$
By the case $p=2$ of Lemma 5, $G$ contains a triangle $C_3.$ Then successively applying Theorem 2 we deduce that $G$ is pancyclic. $\hfill\Box$

{\bf Remark.} The Chv\'{a}tal-Erd\H{o}s theorem on hamiltonian graphs ([4], [1] and [8]) states that for a graph $G,$ if $\kappa(G)\ge \alpha(G)$ then $G$ is hamiltonian,
where $\kappa$ and $\alpha$ denote the connectivity and independence number, respectively. Bondy [2] proved that if a graph satisfies  Ore's condition, then it satisfies the Chv\'{a}tal-Erd\H{o}s condition. A computer search for graphs of lower orders shows that there are many graphs which satisfy the condition in Theorem 6, but do not satisfy the Chv\'{a}tal-Erd\H{o}s condition. There are exactly $398$ such graphs of order $9.$ For every integer $n\ge 7$ we give an example. Let $G_1=K_{n-3}^{-}$ be the graph obtained from $K_{n-3}$ by deleting one edge $xy$ and let $G_2=uvw$ be a triangle that is vertex-disjoint from $G_1.$ Construct a graph $Z_n$ form $G_1$ and $G_2$ by adding two edges $xu$ and $yv.$ The graph $Z_9$ is depicted in Figure 3.
\begin{figure}[h]
\centering
\includegraphics[width=0.3\textwidth]{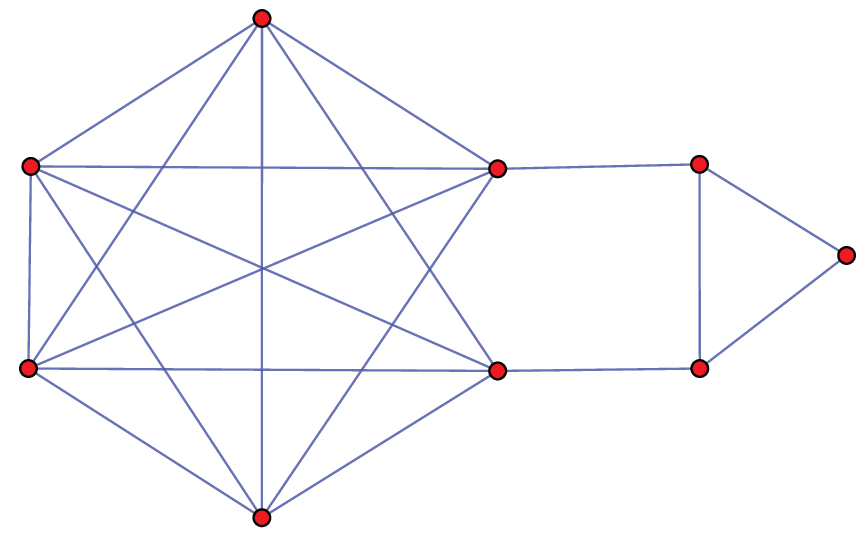}
\caption{The graph $Z_9$}
\end{figure}

Clearly $Z_n$ is a $2$-connected $[4,2]$-graph of order $n,$ but $2=\kappa(Z_n)<\alpha(Z_n)=3.$

\vskip 5mm
{\bf Acknowledgement.} This research  was supported by the NSFC grant 12271170 and Science and Technology Commission of Shanghai Municipality
 grant 22DZ2229014.


\begin{thebibliography}{99}
\bibitem{1} J.A. Bondy, Pancyclic graphs. I, J. Combinatorial Theory Ser. B, 11(1971), 80-84.
\bibitem{2} J.A. Bondy, A remark on two sufficient conditions for Hamilton cycles, Discrete Math., 22(1978), no.2, 191-193.
\bibitem{3} J.A. Bondy and U.S.R. Murty, Graph Theory, GTM 244, Springer, 2008.
\bibitem{4} V. Chv\'{a}tal and P. Erd\H{o}s, A note on Hamiltonian circuits, Discrete Math., 2(1972), 111-113.
\bibitem{5} J.C. George, A. Khodkar and W.D. Wallis, Pancyclic and Bipancyclic Graphs, Springer, 2016.
\bibitem{6} C. Liu and J. Wang, $[s,t]$-graphs and their hamiltonicity, (Chinese), J. Shandong Normal Univ. Nat. Sci., 20(2005), no.1, 6-7.
\bibitem{7} X. Liu, J. Wang and G. Gao, Cycles in $2$-connected $[4,2]$-graphs, (Chinese), J. Shandong Univ. Nat. Sci., 42(2007), no.4, 32-35.
\bibitem{8} D.B. West, Introduction to Graph Theory, Prentice Hall, Inc., 1996.
\bibitem{9} X. Zhan, Extremal numbers of positive entries of imprimitive nonnegative matrices, Linear Algebra Appl. 424(2007), no.1, 132–138.
\end{thebibliography}
\end{document}